\documentclass[runningheads]{llncs}
\usepackage{physics}
\usepackage{enumitem} 
\usepackage{bm}
\usepackage{stmaryrd}
\usepackage[T1]{fontenc}
\usepackage{esvect}
\usepackage{graphicx}
\usepackage{hyperref}
\hypersetup{
   colorlinks   =  true,
    linkcolor    = blue,
    citecolor    = red,
     urlcolor	=blue,     
	urlbordercolor = {1 1 0}
}
\usepackage{color}

\usepackage{amsmath, amsfonts, amssymb}
\renewcommand{\d}{\mathrm{d}}
\begin{document}

\title{Applications of standard and Hamiltonian stochastic Lie systems}
%
%
\author{Javier de Lucas\inst{1}\orcidID{0000-0001-8643-144X} and Marcin Zaj\c ac\inst{1,2}\orcidID{0000-0003-4914-3531}}
\authorrunning{J. de Lucas and  M. Zaj\k ac}
%
\institute{Department of Mathematical Methods in Physics, University of Warsaw, Poland \and
Instituto de Matem\'atica Interdisciplinar, UCM, Madrid,  Spain\\
\email{ }}
\maketitle              
\begin{abstract}
A stochastic Lie  system on a manifold $M$ is a stochastic differential equation whose dynamics is described by a linear combination with functions depending on $\mathbb{R}^\ell$-valued semi-martigales of vector fields on $M$ spanning a finite-dimensional Lie algebra. We analyse new examples of stochastic Lie systems and Hamiltonian stochastic Lie systems, and review and extend the coalgebra method for Hamiltonian stochastic Lie systems. We apply the theory to  biological and epidemiological models, stochastic oscillators, stochastic Riccati equations, coronavirus models, stochastic Ermakov systems, etc.
\end{abstract}


\section{Introduction}
A {\it Lie system} is a $t$-dependent system of ordinary differential equations in normal form whose general solution can be written as a $t$-independent function, the {\it superposition rule}, of a generic family of particular solutions and some constants related to initial conditions \cite{CGM00,CGM07,Lie1893,LS21}. Superposition rules occur in approximate and numerical analysis \cite{LS21,PW83}.  The Lie--Scheffers theorem states that a $t$-dependent system of ordinary differential equations on an $n$-dimensional manifold $M$ in normal form, 
$   \frac{\text{d} \Gamma}{\text{d} t} = X(t,\Gamma ), 
$ which amounts to a  $t$-dependent vector field $X=\sum_{i=1}^n X^i(t,\Gamma)\partial/\partial \Gamma^i$ on $M$, 
admits a superposition rule if and only if $X=\sum_{\alpha=1}^rb_\alpha(t)X_\alpha$ for vector fields $X_1,\ldots,X_r$ on $M$ spanning an $r$-dimensional Lie algebra of vector fields and certain functions $b_1(t),\ldots,b_r(t)$. Lie systems are not commonly found among differential equations, but their significant applications and relevant mathematical features strongly justify their study \cite{CFH_23,CGM00,LS21,PW83}.

Efforts have been made to extend the theory of Lie systems and superposition rules to more general settings, including $t$-dependent Schr\"odinger equations \cite{Dissertationes}, partial differential equations \cite{AFV_09,CGM07,OG_00}, quasi--Lie systems and families \cite{Dissertationes}, et cetera \cite{CGM00}.  The theory of Lie systems and their superposition rules was extended to stochastic differential equations (SDEs) in \cite{LO09} and reviewed and extended in \cite{LSR_25}, where technical details were fixed, Hamiltonian stochastic (HS) Lie systems were defined, and a Poisson coalgebra method for some of them, commented. 

This short work shows relevant new examples and potential applications of stochastic Lie systems, in general, and HS Lie ones, in particular. These include epidemic, biological, physical, economical, and coronavirus models.  This enlarges the applications in \cite{FLRZ_25,LO09} given by geometric Brownian motions, inhomogeneous systems, stochastic oscillators, epidemic models, and stochastic Wei--Norman equations. Stochastic Lie systems together with their Hamiltonian versions are applicable to everyday problems, like Black-Scholes pricing \cite{FLRZ_25,LO09} or cases of coronavirus models, as shown here for the first time. SDEs complement deterministic equations \cite{CFH_23} describing phenomena beyond the reach of deterministic models, like outbreak likelihood  \cite{Al09}. Other potential applications of our methods to coronavirus models, biological systems, and other topics are briefly discussed \cite{Ar03,BJBA_21}. This work also   slightly generalises the stochastic Poisson coalgebra method to cover general HS Lie systems, namely with coefficients depending on stochastic variables and the time instead of only the time as in \cite{FLRZ_25}. Then, we apply it to stochastic Riccati equations, diffusion Lotka--Volterra models \cite{Ar03}, and  coronavirus models to derive new superposition rules. 

Let us introduce HS Lie systems \cite{Ar74,FLRZ_25,LO08}. Consider the SDE on $M$ given by 
\begin{equation}\label{Eq:practicalSDE}
\delta \Gamma^i=X^i_1(B,\Gamma)\delta t+\sum_{\alpha=2}^\ell X^i_\alpha(B,\Gamma)\circ \delta B^\alpha\,,\qquad i=1,\ldots,n\,,
\end{equation}
where $X^i_1 ,\ldots,X_\ell^i \in C^\infty(\mathbb{R}^\ell\times M)$,  $B=(B^1=t,\ldots,B^\ell)$, $\Gamma=(\Gamma^1,\ldots,\Gamma^n)\in M$, and $i=1,\ldots,n$. The symbol $\circ$ indicates that \eqref{Eq:practicalSDE} is a Stratonovich SDE. 
More precisely, $(\Omega, \mathcal{F}, P)$ is a {\it probability space}, where $\Omega$ is a manifold, $\mathcal{F}$ is a $\sigma$-algebra of subsets of $\Omega$, the $P:\mathcal{F}\rightarrow [0,1]$ is a probability function,  and  $B:\mathbb{R}_+\times\Omega\rightarrow \mathbb{R}^\ell$ is a {\it semi-martingale}. Brownian motions,  Wiener processes, and other common stochastic processes are examples of semi-martingales. Applications of SDEs are frequently written via It\^o calculus \cite{Al09,WCDG18}. If  $X^i_\alpha=X^i_\alpha(t,\Gamma)$ for $\alpha=2,\ldots,\ell$ and $i=1,\ldots,n$, the Stratonovich SDE \eqref{Eq:practicalSDE} amounts to an It\^o SDE
\begin{equation}\label{Eq:TranItoStratonovich}
\delta \Gamma^i=\left(X^i_1(B,\Gamma)-\frac 12\sum_{\beta=2}^\ell\sum_{j=1}^n\frac{\partial X_\beta^i}{\partial \Gamma^j}(t,\Gamma)X_\beta^j(t,\Gamma)\right)\delta t +\sum_{\beta=2}^\ell X^i_\beta(t,\Gamma )\delta {B}^\beta_t.
\end{equation} A more general dependence of the $X^i_\alpha$ gives new terms above. 
Rewrite \eqref{Eq:practicalSDE}  as  
\begin{equation}
  \label{stochastic differential equation expression} \delta \Gamma = \mathfrak{S} (B,
  \Gamma) \circ\delta B\,,
\end{equation}
where $\mathfrak{S} (B, \Gamma) : T_B \mathbb{R}^\ell
\to T_\Gamma M$, with $(B,\Gamma)\in \mathbb{R}^\ell\times M$, describes a Stratonovich operator. A basis in
$T^{\ast} \mathbb{R}^\ell$ identifies $\mathfrak{S}(B,
\Gamma)$ with its $\ell$-tuple of components $(\mathfrak{S}_1 (B, \Gamma), \ldots, \mathfrak{S}_\ell (B, \Gamma))$ in the chosen basis. Every particular solution to \eqref{stochastic differential equation expression} is also a semi-martingale $
\Gamma:\mathbb{R}_+\times \Omega\rightarrow M$. A particular solution has initial condition $\Gamma_0
\in M$ when $\Gamma(0,\omega_0)=\Gamma_0$ for every $\omega_0\in \Omega$ with probability one. Time is considered as the first component $B_1:(t,\omega_0)\in \mathbb{R}\times \Omega\mapsto t\in\mathbb{R}$ of $B$.  

The solution to \eqref{stochastic differential equation expression} is described via {\it Stratonovich integrals} as
\begin{equation}\label{Eq:SolStr}
\Gamma(t)-\Gamma(0)=\int_0^t\mathfrak{S}_1(B,\Gamma)\delta t+\sum_{\beta=2}^{\ell}
\int_0^t\mathfrak{S}_\beta(B,\Gamma) \circ \delta B^\beta_t.
\end{equation}
The form of \eqref{Eq:TranItoStratonovich} is due to the change from It\^o to Stratonovich integrals in \eqref{Eq:SolStr}.

\begin{definition}\label{Def:StoLieSys} A {\it stochastic Lie system} on $M$ is a SDE \eqref{stochastic differential equation expression}  so that $B:\mathbb{R}_+\times \Omega\rightarrow \mathbb{R}^\ell$ is a semi-martingale and $\mathfrak{S}$ is a Stratonovich operator 
\begin{equation}
\label{Eq:STL}\mathfrak{S}(B,\Gamma)=\left(\sum_{\alpha=1}^rb_1^\alpha(B)Y_\alpha(\Gamma),\ldots,\sum_{\alpha=1}^rb_\ell^\alpha(B)Y_\alpha(\Gamma)\right),\,\, \forall\Gamma\in M,\,\, \forall B\in \mathbb{R}^\ell,
\end{equation}
for $b^a_\alpha:B\in \mathbb{R}^\ell\mapsto b^a_\alpha(B)\in \mathbb{R}$, with $\alpha=1,\ldots,r,$ $a=1,\ldots,\ell$, and an $r$-dimensional Lie algebra of vector fields $\langle Y_1,\ldots,Y_r\rangle $, called a {\it Vessiot--Guldberg (VG) Lie algebra} of the stochastic Lie system given by \eqref{Eq:STL}.
\end{definition}

Note that \eqref{Eq:STL} can be considered as an $\ell$-element family of $\mathbb{R}^\ell$-dependent vector fields on $M$. 
Meanwhile, 
HS Lie systems are as follows.

\begin{definition} A {\it Hamiltonian stochastic (HS) Lie system} on $M$ is a stochastic Lie system $\delta\Gamma=\mathfrak{S}(B,\Gamma)\circ \delta B$ admitting a VG Lie algebra $V$ of Hamiltonian vector fields relative to a geometric structure on $M$. A Lie algebra containing Hamiltonian functions for all elements of $V$ is called a {\it Lie--Hamilton Lie algebra}.
\end{definition}


\section{Superposition rules and Poisson coalgebra method}\label{Sec:SSRandLT}

Let us introduce superposition rules  for SDEs and our slight extension of the Poisson coalgebra method to general SH Lie systems, i.e. for $\mathfrak{S}(B,\Gamma)$ instead of only $\mathfrak{S}(t,\Gamma)$ \cite{FLRZ_25,LO09}.  We write $\mathfrak{X}(M)$ for the space of vector fields on $M$.

\begin{definition}\label{def 1}
    A {\it superposition rule} for \eqref{stochastic differential equation expression} is a function $\Phi: M^{m + 1} \to M$ such that the general solution $\Gamma$ to~\eqref{stochastic differential equation expression} reads
    $ \Gamma = \Phi (z
    ; \Gamma_{(1)}, \ldots, \Gamma_{(m)})$, for all $z\in M$ and generic particular solutions $\Gamma_{(1)}, \ldots, \Gamma_{(m)}:\mathbb{R}\times \Omega\rightarrow M$ of~\eqref{stochastic differential equation expression}, 

\end{definition}
A superposition rule $\Phi$ does not explicitly depend on either $B\in \mathbb{R}^\ell$ or $\Omega$. Let us briefly describe how to derive superposition rules for stochastic Lie systems.

The {\it diagonal prolongation to $M^k$ of a vector bundle} $\tau:F\to M$  is the vector bundle $\tau^{[k]}:(f_{(1)},\dotsc,f_{(k)})\in F^k \mapsto (\tau(f_{(1)}),\dotsc,\tau(f_{(k)}))\in M^k$. Each section $e:M\to F$ of  $\tau$ has a  {\it  diagonal prolongation} to a section $e^{[k]}$ of $\tau^{[k]}$ given by
\begin{equation*}
    e^{[k]}(\Gamma_{(1)},\ldots,\Gamma_{(k)})=(e(\Gamma_{(1)}),\ldots ,e(\Gamma_{(k)}))\,,\qquad \forall (\Gamma_{(1)},\ldots,\Gamma_{(k)})\in M^k.
\end{equation*}
If $e(\Gamma_{(a)})$ is assumed to take values in the $a$-th copy of $F$ within $F^k$, one can simply write
$ e^{[k]}(\Gamma_{(1)},\ldots,\Gamma_{(k)})=\sum_{a=1}^ke(\Gamma_{(a)})\,$. The {\it  diagonal prolongation of  $f\in C^\infty (M)$} to $M^k$ is 
$f^{[k]}:(\Gamma_{(1)},\ldots,\Gamma_{(k)})\in M^k\mapsto  \sum_{a=1}^kf(\Gamma_{(a)})\in \mathbb{R}\,.
$ The canonical isomorphisms
$(T M)^{[k]}\simeq T M^k$  and $ (T^* M)^{[k]}\simeq T^* M^k\,$ allow us to understand
 $X^{[k]}$ for $X\in\mathfrak{X}(M)$ as a vector field on $M^k$, and the diagonal prolongation, $\alpha^{[k]}$, of a one-form $\alpha$ on $M$ as a one-form $ {\alpha}^{[k]}$ on $M^k$. In fact,  if $Y\in \mathfrak{X}(M)$, then   
  $
        {Y}^{[k]}(\Gamma_{(1)},\ldots,\Gamma_{(k)})=\sum_{a=1}^kY(\Gamma_{(a)})\in \mathfrak{X}(M^k)$. In particular, $Y=\Gamma\partial/\partial \Gamma\in \mathfrak{X}(\mathbb{R})$ gives $Y^{[k]}=\sum_{a=1}^k\Gamma_{(a)}\partial/\partial \Gamma_{(a)}\in \mathfrak{X}(\mathbb{R}^k)$. Note that $X\in \mathfrak{X}(M)\mapsto {X}^{[k]}\in \mathfrak{X}(M^{k})$ is a Lie algebra morphism.

The diagonal prolongation of an  $\mathbb{R}^\ell$-dependent vector fields on $M$, namely a mapping $X:\mathbb{R}^\ell\times M\rightarrow T M$ such that $X_B=X(B,\cdot)$ is a standard vector field on $M$ for every $B\in \mathbb{R}^\ell$,  is the $\mathbb{R}^\ell$-dependent vector field ${X}^{[k]}$ on $M^k$ whose value for every fixed $B\in \mathbb{R}^\ell$, let us say ${X}^{[k]}_B$, is  $(X_B)^{[k]}$.

 The proof of the stochastic Lie theorem shows that a superposition rule for stochastic Lie systems can be obtained almost like for Lie systems (cf. \cite{CGM07,FLRZ_25}):
	\begin{enumerate}
		\item Let $V=\langle X_1,\ldots,X_r\rangle$ be an $r$-dimensional VG Lie algebra for \eqref{stochastic differential equation expression}.
  \item Find the smallest $m\in \mathbb{N}$ so that $X^{[m]}_1\wedge \ldots\wedge X_r^{[m]}\neq 0$  
		at a generic point.
		\item Set coordinates $\Gamma^1,\ldots,\Gamma^n$ on $M$ and  use them in each copy of $M$ within $M^{m+1}$ to get a coordinate system
		$\{\Gamma^i_{(a)}\mid i=1,\ldots,n,\,\,a=0,\ldots,m\}$ on $M^{m+1}$.
		Obtain first integrals $F_1,\ldots, F_n$ for  all
		$X^{[m+1]}_1,\ldots,  X^{[m+1]}_r$ so that 
		\begin{equation}\label{Cond}
		\frac{\partial(F_1,\ldots,F_n)}{\partial(\Gamma_{(0)}^1,\ldots,\Gamma_{(0)}^n)}\neq 0.
		\end{equation}
		  
		\item 
	The equations $F_i=k_i$, for $i=1,\ldots,n$, enable us to express  $\Gamma_{(0)}^1\ldots,\Gamma_{(0)}^n$ in
		terms of remaining $\Gamma^1_{(a)},\ldots,\Gamma_{(a)}^n$ and $k_1,\ldots,k_n$, giving a superposition rule for \eqref{Eq:STL} depending on $m$ particular solutions
		and $k_1,\ldots, k_n$.  
	\end{enumerate}

Let us review and slightly generalise the Poisson coalgebra method for HS Lie systems devised in \cite{FLRZ_25}. 
Every HS Lie system is related to a Stratonovich operator $\mathfrak{H}$ given by $\ell$ components and each one can be understood as an $\mathbb{R}^\ell$-dependent vector field. Hence, $\mathfrak{H}$, for every $B\in \mathbb{R}^\ell$, is a family of $\ell$ vector fields on $M$. In this manner, one can prolong it to a section $\mathfrak{H}^{[m]}$ of $\mathbb{R}^\ell\times \oplus^\ell (TM^{[m]})\rightarrow \mathbb{R}^\ell\times M^{[m]}$. 
Moreover, $\mathfrak{H}$ is Hamiltonian, and let us assume that it is so relative to a symplectic form. Then, one obtains a mapping ${h}:(B,\Gamma)\in \mathbb{R}^\ell\times M\mapsto (h_1(B,\Gamma),\ldots,h_\ell(B,\Gamma))\in \mathbb{R}^\ell$. Then, $h$ can be extended to a mapping  $h^{[m]}:\mathbb{R}^\ell \times M^m\rightarrow \mathbb{R}^\ell$ so that  $h_\alpha^{[m]}(B,\cdot)$ is the extension of the  Hamiltonian function $h_\alpha(B,\cdot)$ for every $B\in \mathbb{R}^\ell$ and $\alpha=1,\ldots,\ell$.  Then, the following result is immediate \cite{FLRZ_25,LS21}.

\begin{proposition}\label{Cor:LieAlgJL2}
    If $\mathfrak{H}$ is a HS Lie system with an $\mathbb{R}^\ell$-dependent Hamiltonian $\widetilde{h}:\mathbb{R}^\ell \times M\rightarrow \mathbb{R}^\ell$ relative to a symplectic form $\omega$ on $M$, then $\mathfrak{H}^{[m]}$ is a HS Lie system relative to $\omega^{[m]}$ with a Hamiltonian $\widetilde{h}^{[m]}:\mathbb{R}^\ell\times M^m\rightarrow \mathbb{R}^\ell$. If $h_1,\ldots,h_r$ is a basis of a Lie--Hamilton Lie algebra $\mathfrak{W}$ for $\mathfrak{H}$, then $h_1^{[m]},\ldots,h_r^{[m]}$ is a basis of the Lie--Hamilton Lie algebra for $\mathfrak{H}^{[m]}$.  Then, $f\in C^\infty(M)$ is a constant of the motion for $\mathfrak{H}$ if Poisson commutes with the elements of $\{h_t\}_{t\in \mathbb{R}}$.  Let $\{v_1,\ldots ,v_r\}$ be a basis of linear coordinates on $\mathfrak{g}^*\simeq \mathfrak{W}$.  
    If $C$ is a Casimir function  on $\mathfrak{g}^\ast$ and
         $C=C(v_1,\ldots,v_r)$, then the following are  constants of motion of $\mathfrak{H}^{[m]}$: 
    \begin{equation}\label{Eq:ExCas}
C\left(\sum_{a=1}^sh_1(x_{(a)}),\ldots, \sum_{a=1}^sh_r(x_{(a)})\right)\,,\qquad 1\leq s\leq m.
\end{equation}
\end{proposition}
\section{New applications of stochastic Lie systems}

Let us consider the stochastic differential equation
\begin{equation}\label{Eq:StoRicc}
\delta \Gamma=(\Gamma^2+f\Gamma+g)\delta t,\qquad \Gamma\in \mathbb{R},
\end{equation}
where $f,g$ are  well-behaved stochastic processes depending on $\Gamma$ and the time $t$, e.g. functions of semi-martingales \cite{BA_80,BSV_72}. Note that the form of \eqref{Eq:StoRicc} in the It\^o and the Stratonovich forms is the same.  Equations \eqref{Eq:StoRicc} are mathematically and physically interesting, e.g. they occur in the study of stochastic harmonic oscillators \cite{Ar03,BA_80}. In particular, during the COVID lockdown period, deterministic Riccati equations were utilized for its analysis
\cite{FDK_21}. Moreover, \eqref{Eq:StoRicc} is a particular case of the stochastic matrix Riccati equations with stochastic coefficients \cite{BSV_72,BMN_20}. 

More general stochastic Riccati equations are given by 
\begin{equation}\label{eq:StoRicc}
\delta \Gamma=(b_2(B)\Gamma^2+b_1(B)\Gamma+b_0(B))\delta t+(b'_2(B)\Gamma^2+b'_1(B)\Gamma+b'_0(B))\circ \delta \mathcal{B} 
\end{equation}
for a Brownian motion $\mathcal{B}$, functions $b_0,b_1,b_2,b'_0,b'_1,b'_2\in C^\infty(\mathbb{R}^2)$ for $B=(t,\mathcal{B})$, and $\Gamma\in \mathbb{R}$. Systems \eqref{eq:StoRicc}  retrieve as particular cases affine stochastic differential equations appearing in SIR models (the equation for the $R$ variable),  some geometric Brownian motions \cite{LO08}, etc.  

It is worth considering a stochastic harmonic oscillator of the form 
\begin{equation}\label{eq:SHO}
\delta x=v\delta t,\quad 
\delta v=-(\omega^2(t,\mathcal{B})x+\gamma(t,\mathcal{B})v)\delta t-(\omega_B^2(t,\mathcal{B})x+\gamma_B(t,\mathcal{B})v)\circ \delta \mathcal{B}
\end{equation}
for arbitrary functions $\gamma,\gamma_B,\omega,\omega_B:\mathbb{R}^2\rightarrow \mathbb{R}$. 
Introducing $\Gamma=x/v$, one has that
$$
\delta \Gamma=(1+\gamma(t,\mathcal{B})\Gamma+\omega^2(t,\mathcal{B})\Gamma^2)\delta t+(\gamma_B(t,\mathcal{B})\Gamma+\omega^2_B(t,\mathcal{B})\Gamma^2)\circ \delta B.
$$
The stochastic harmonic oscillator presented above is one of the many instances of stochastic harmonic oscillators currently analysed \cite{CEGR_14,FSZXL_10}. This also shows the relevance of \eqref{eq:StoRicc}. It is worth noting, as this seems to be absent in the present literature, that the addition of a drift in \eqref{eq:SHO} is incompatible with obtaining a stochastic Riccati equation for $\Gamma$.

SDEs \eqref{eq:StoRicc} are stochastic Lie systems related to a VG Lie algebra spanned by $X_\alpha=\Gamma^\alpha \frac{\partial}{\partial \Gamma}$ for $\alpha=0,1,2$, which isomorphic to $\mathfrak{sl}_2$. In fact,
$$
[X_0,X_1]=X_0,\qquad [X_0,X_2]=2X_1,\qquad [X_1,X_2]=X_2.
$$
Let us apply our Poisson coalgebra method to it. Then,  $X^{[2]}_0\wedge X^{[2]}_1\wedge X^{[2]}_2=0$ almost everywhere, but $X_0^{[3]}\wedge X_1^{[3]}\wedge X_2^{[3]}\neq 0$ at a generic point on $\mathbb{R}^3$. Then, a superposition rule with three particular solutions exists and, to derive it, consider  
$
X^{[4]}_\alpha =\sum_{a=0}^3\Gamma_{(a)}^\alpha\frac{\partial}{\partial\Gamma_{(a)} },$ for $\alpha=0,1,2$, 
which are Hamiltonian relative to the symplectic form $\omega=\sum_{i=0}^1 \frac{d\Gamma_{(2i)}\wedge d\Gamma_{(2i+1)}}{\Gamma_{(2i)}-\Gamma_{(2i+1)}}$ with Hamiltonian functions 
$$
\widetilde{h}_0=\sum_{i=0}^1\frac{1}{\Gamma_{(2i)}-\Gamma_{(2i+1)}},\,
 \widetilde{h}_1=\frac 12\sum_{i=0}^1\frac{\Gamma_{(2i)}+\Gamma_{(2i+1)}}{\Gamma_{(2i)}-\Gamma_{(2i+1)}},\,\,\widetilde{h}_2=\sum_{i=0}^1\frac{\Gamma_{(2i)}\Gamma_{(2i+1)}}{\Gamma_{(2i)}-\Gamma_{(2i+1)}}.
$$
Note that
$$
\{\widetilde{h}_0,\widetilde{h}_1\}=-\widetilde{h}_0,\qquad \{\widetilde{h}_0,\widetilde{h}_2\}=-2\widetilde{h}_1,\qquad \{\widetilde{h}_1,\widetilde{h}_2\}=-\widetilde{h}_2.
$$
They are indeed diagonal prolongations to $\mathbb{R}^4$ from the Hamiltonian functions of the prolongation to $\mathbb{R}^2$ of \eqref{eq:StoRicc} and a symplectic form, which are quite straightforward. 
Then, $\{\widetilde{h}_i,\mathcal{C}\}=0$ for $i=0,1,2$ and $\mathcal{C}=\widetilde{h}_2\widetilde{h}_0-\widetilde{h}_1^2$, and $\mathcal{C}$ becomes a first integral for $X^{[4]}_0,X^{[4]}_1,X^{[4]}_2$  satisfying \eqref{Cond} and gives a superposition rule for  \eqref{eq:StoRicc} given by
$
\Phi_{\rm Ric}:\mathbb{R}^3\times \mathbb{R}\rightarrow \mathbb{R}
$
of the form 
\begin{equation}
    \Phi_{\rm Ric}(\Gamma_{(1)},\Gamma_{(2)},\Gamma_{(3)},z)=\frac{\Gamma_{(3)}(\Gamma_{(1)}-\Gamma_{(2)}) + z\Gamma_{(1)}(\Gamma_{(3)} -\Gamma_{(2)} )}{\Gamma_{(1)} -\Gamma_{(2)}  +z(\Gamma_{(3)} -\Gamma_{(2)} )}\,, \label{solutin1}
\end{equation}
which implies that the general solution for \eqref{Eq:StoRicc}, let us say $\Gamma(t)$, can be brought into the form $\Gamma(t)=\Phi_{\rm Ric}(\Gamma_{(1)}(t),\Gamma_{(2)}(t),\Gamma_{(3)}(t),z)$, where 
 $\Gamma_{(1)}(t),\Gamma_{(2)}(t),\Gamma_{(3)}(t)$ are three different particular solutions of \eqref{eq:StoRicc} and $z\in \mathbb{R}$. This expression is similar to the known superposition rule for Riccati equations, but applies for stochastic Riccati equations, which are more general. This is the first time a superposition rule has been derived for stochastic Riccati equations. 

Let us consider now a stochastic Ermakov system \cite{CEGR_14} of the form
\begin{equation}\label{eq:ErmakovSto}
\delta \rho=v\delta t,\qquad \delta v=\left(-\omega^2(t,\mathcal{B})\rho+\frac{k}{\rho^3}\right)\delta t+\sigma \rho \circ \delta \mathcal{B},
\end{equation}
for a certain function $\omega\in C^\infty(\mathbb{R}^2)$ and  constants $\sigma,k\in \mathbb{R}$. 
Physically, this is an isotropic oscillator on $\mathbb{R}^3$ with a perturbation stochastic term. 

System \eqref{eq:ErmakovSto} is a stochastic Lie system associated with a Vessiot--Guldberg Lie algebra isomorphic to $\mathfrak{sl}_2$ spanned by the basis
$$
X_1=-\rho\frac{\partial}{\partial v},\qquad X_2= \frac 12\left(v\frac{\partial}{\partial v}-\rho\frac{\partial}{\partial \rho}\right),\qquad X_3=v\frac{\partial}{\partial \rho}+\frac{k}{\rho^3}\frac{\partial}{\partial v},
$$
with commutation relations
$$
[X_1,X_2]=X_1,\qquad [X_1,X_3]=2X_2,\qquad [X_2,X_3]=X_3.
$$
Then $\langle X_1,X_2,X_3\rangle$ is a Lie algebra of Hamiltonian vector fields relative to the symplectic form $\omega=d\rho\wedge d v$. In fact, $X_1,X_2,X_3$ have  Hamiltonian functions, respectively, given by
$$
h_1=\frac 12 \rho^2,\qquad h_2=-\frac 12\rho v,\qquad h_3=\frac 12\left(v^2+\frac{k}{\rho^2}\right).
$$
They close a Lie algebra isomorphic to $\mathfrak{sl}_2$. It follows that one may have a Casimir of $\mathfrak{sl}_2$, which gives rise to $\mathcal{C}=h_1h_3-h_2^2$.
As in the classical Ermakov system with no term $\delta \mathcal{B}$, one obtains that the SH Lie system \eqref{eq:ErmakovSto} admits a superposition rule depending on two particular solutions generated by the extension to two three copies of the system. The functions induced by a Casimir of $\mathfrak{sl}_2$ via the Poisson coalgebra method (with our correction) given by \cite{CL_08} 
$$F_1=\left(\sum_{a=0}^1h_1(\rho_{(a)},v_{(a)})\right)\left(\sum_{a=0}^1h_3(\rho_{(a)},v_{(a)})\right)-\left(\sum_{a=0}^1 h_2(\rho_{(a)},v_{(a)})\right)^2 ,$$
$$F_2=\left(\sum_{a=0,2}h_1(\rho_{(a)},v_{(a)})\right)\left(\sum_{a=0,2}h_3(\rho_{(a)},v_{(a)})\right)-\left(\sum_{a=0,2}h_2(\rho_{(a)},v_{(a)})\right)^2.
$$
 Affine models on $\mathbb{R}^7$ with three stochastic Brownian variables $\mathcal{B}_1,\mathcal{B}_2,\mathcal{B}_3$ related to coronavirus models can be found in \cite{BJBA_21}. For certain limit cases of the parameters, some of the variables are given by the stochastic Lie system on $\mathbb{R}^2_+$ of the form
\begin{equation}\label{eq:Model1}
    \delta H=-A(t)H\delta t- B(t)H\delta \mathcal{B}_1,\qquad \delta R=-A(t)R\delta t+B(t)H\delta \mathcal{B}_1,
\end{equation}
for certain $t$-dependent functions $A(t),B(t)$, 
which is Hamiltonian with a non-abelian two-dimensional VG Lie algebra  $\langle M_1=H(\partial_H-\partial_R),M_2=-H\partial_H-R\partial_R\rangle$ of Hamiltonian vector fields relative to $\omega=dH\wedge dR/(RH+H^2)$ and Hamiltonian functions $h_1=\ln(H+R)$ and $h_2=\ln( H/(H+R))$. Note that
$$
\{h_1,h_2\}=M_2h_1=-1.
$$
The system \eqref{eq:Model1} admits a superposition rule depending on one particular solution. The functions $h^{[2]}_2=\ln(H_{(0)}/(H_{(0)}+R_{(0)}))+\ln(H_{(1)}/(H_{(1)}+R_{(1)}))$ and $h^{[2]}_1=\ln(H_{(0)}+R_{(0)})+\ln(H_{(1)}+R_{(1)}),
$ are the Hamiltonian functions of  $M_1^{[2]},M_2^{[2]}$. Since $\{h_\alpha^{[2]},h_\beta(H_0,R_0)-h_\beta(H_1,R_1)\}=0$ for $\alpha,\beta=1,2$,  there are common first integrals for $M_1^{[2]},M_2^{[2]}$ of given by
$$F_1=\ln(H_{(0)}+R_{(0)})-\ln(H_{(1)}+R_{(1)}),
$$
$$
F_2=\ln(H_{(0)}/(H_{(0)}+R_{(0)}))-\ln(H_{(1)}/(H_{(1)}+R_{(1)})),
$$
which satisfy \eqref{Cond} and give rise to a superposition rule from $F_1=\ln k_1$ and $F_2=\ln k_2$ depending on two parameters $k_1,k_2\in \mathbb{R}$ of the form
$$
H_0 = k_1 k_2 H_1,\qquad
R_0 =  k_1(H_1 + R_1) - k_1 k_2 H_1.
$$

Let us now consider the stochastic Lotka--Volterra system with diffusion \cite{Ar03}
\begin{equation}\label{Eq:N1N2}
\begin{gathered}
\delta N_1  = (b_1-a_1N_2)N_1\delta t+\sigma_1 N_1\delta \omega_1\,,\qquad 
\delta N_2  = b_2N_2\delta t+\sigma_2N_2\delta \omega_2\,,
\end{gathered}
\end{equation}
on $\mathbb{R}_+^2$, where $b_1,a_1,\sigma_1$ are constants and $\omega_1,\omega_2$ are Brownian motions,
that one can find as a particular case of the system analysed in \cite{Ar03} and studied in \cite{FLRZ_25} from the perspective of HS Lie systems. System (\ref{Eq:N1N2}), even when $b_1=b_1(t),a_1=a_1(t),\sigma_1=\sigma_1(t)$, becomes a stochastic Lie system related to a VG Lie algebra, spanned by
$
Z_1=N_1\frac{\partial}{\partial N_1},\,Z_2=N_2\frac{\partial}{\partial N_2}\,,Z_3=N_1N_2\frac{\partial}{\partial N_1}\,,
$
isomorphic  to $\mathbb{R}^2\ltimes \mathbb{R}$ consisting of Hamiltonian vector fields relative to 
$
\omega=\frac{\d N_1\wedge \d N_2}{N_1N_2}.
$ A superposition rule can be found for $m=2$ and the superposition rule follows from  obtaining two first integrals for $Z_1^{[3]},Z_2^{[3]},Z_3^{[3]}$ satisfying \eqref{Cond}. Our previous results show  that the diagonal prolongation of \eqref{Eq:N1N2} to $(\mathbb{R}_+^2)^2$ is Hamiltonian  relative to $\omega^{[2]}$, and it admits a Hamiltonian Lie symmetry 
$$
\frac{(N_2)_{(0)}}{(N_2)_{(1)}}\left((N_1)_{(0)}\frac{\partial}{\partial (N_1)_{(0)}}-(N_1)_{(1)}\frac{\partial}{\partial (N_1)_{(1)}}\right),
$$
with a Hamiltonian function $F_1=(N_2)_{(0)}/(N_2)_{(1)}$, which becomes  a first integral of $Z^{[2]}_1,Z^{[2]}_2,Z^{[2]}_3$ and, hereafter, of $Z^{[3]}_1,Z^{[3]}_2,Z^{[3]}_3$. To obtain a second first integral for $Z^{[3]}_1,Z^{[3]}_2,Z^{[3]}_3$ satisfying \eqref{Cond}, one has to consider that as diagonal prolongations are invariant relative to interchanges $\Gamma_{(a)}\leftrightarrow\Gamma_{(b)}$, this transformation maps first integrals of diagonal prolongations into first integrals of diagonal prolongations and $x=(N_2)_{(1)}/(N_2)_{(2)}$ is also a first integral of $Z^{[3]}_1,Z^{[3]}_2,Z^{[3]}_3$. Using the previous first integrals with $t=(N_2)_{(1)},u=(N_1)_{(1)},y=(N_1)_{(0)}/(N_1)_{(1)},z=(N_1)_{(1)}/(N_1)_{(2)}$, one may write $Z^{[3]}_3$ in these coordinates to obtain
$$
Z_3^{[3]}=t\left[y(F_1-1)\frac{\partial}{\partial y}+z(1-1/F_2)\frac{\partial}{\partial z}+u\frac{\partial}{\partial u}\right]
$$
and restricting it to common first integrals of $Z_1^{[3]},Z_2^{[3]}$, given by functions of $F_1,F_2,z,t$, we derive, e.g. via the characteristic method, a new first integral for $Z_1^{[3]},Z_2^{[3]},Z_3^{[3]}$ given by 
\[F_2=\frac{\left( (N_{2})_{(1)}-(N_{2})_{(2)}\right) \log\left( \frac{(N_{1})_{(0)}}{(N_{1})_{(1)}} \right) - \left( (N_{2})_{(0)} - (N_{2})_{(1)}\right) \log\left( \frac{(N_{1})_{(1)}}{(N_{1})_{(2)}} \right)}{(N_{2})_{(2)}}
.\]
The equations $F_1=\xi_1\in \mathbb{R}_+$ and $F_2=\ln \xi_2$, with $\xi_2\in \mathbb{R}_+$ allow us to write the superposition rule for the initial system $\Phi:(\mathbb{R}_+^2)^2\times \mathbb{R}^2_+\rightarrow \mathbb{R}_+^2$ of the form 
$$
N_1=\xi_2^{\frac{(N_2)_{(2)}}{{(N_2)_{(1)}}- {(N_2)_{(2)}}}}  {(N_1)_{(1)}} \left(\frac{ {(N_1)_{(1)}}}{ {(N_1)_{(2)}}}\right)^{\frac{( {\xi_1}-1)  {(N_2)_{(1)}}}{ {(N_2)_{(1)}}- {(N_2)_{(2)}}}},\qquad \,N_2=\xi_1 (N_2)_{(1)}.
$$

Finally, consider the stochastic Lotka--Volterra system \cite{Ar03} 
\begin{equation}\label{Eq:N1N2bis}
\begin{gathered}
\delta N_1  = (b_1(t)-a_1(t)N_2)N_1\delta t\,,\qquad 
\delta N_2  = b_2(t)N_2\delta t+\sigma_2(t)\delta \omega_2\,,
\end{gathered}
\end{equation}
for arbitrary $t$-dependent functions $b_1(t),b_2(t),a_1(t),\sigma_2(t)$ and admitting a VG Lie algebra 
$
\left\langle N_1\frac{\partial}{\partial N_1},N_2\frac{\partial}{\partial N_2}\,,N_1N_2\frac{\partial}{\partial N_1}\,,\frac{\partial}{\partial N_2}\,\right\rangle.
$
No symplectic form turns \eqref{Eq:N1N2bis} into a Hamiltonian system \cite{LS21}. 

\section{Conclusions}

Among many other new results, the recent paper \cite{FLRZ_25} introduced HS Lie systems and gave the general scheme of a Poisson coalgebra stochastic method to obtain superposition rules for such systems provided their  Stochastic Stratonovich operator is of the form $\mathfrak{S}(t,\Gamma)$. The Poisson coalgebra method has been here explained in detail and slightly extended to cover also the case where the Stochastic Stratonovich operator is of the form $\mathfrak{S}(B,\Gamma)$, instead of the particular case $\mathfrak{S}(t,\Gamma)$ in \cite{FLRZ_25}. The most important advance of this paper is that it develops many new relevant applications and explicit superposition rules for HS Lie systems via the Poisson coalgebra method for the first time. In particular, models  \eqref{eq:ErmakovSto}, \eqref{eq:Model1}, \eqref{Eq:N1N2} and \eqref{Eq:N1N2bis} are studied from the perspective of stochastic Lie systems for the first time. Relevantly, part of these examples introduce for the first time stochastic Lie systems related to coronavirus models, stochastic Erkmakov systems and other biological and epidemic models. 
Meanwhile, the stochastic Riccati equation  \eqref{eq:StoRicc} is an extension of the one briefly analysed in \cite{FLRZ_25}. In particular, this work shows a more detailed analysis of the stochastic Riccati equation in practical problems and it also calculates its superposition rule via our Poisson coalgebra method for the first time. The relation of our stochastic Riccati equation with certain stochastic oscillators has been analysed.  This work also shows the possibility of extending our techniques to stochastic matrix Riccati equations and explain their potential interest in the literature \cite{BSV_72}. We present many new potential applications of HS Lie systems to be further analysed in the future. It is worth stressing that it is important in the theory of Lie systems and their generalisations to extend the number of examples to new realms, as done in this work.

\section{Acknowledgments}
M. Zaj\k ac acknowledges funding from his NCN PRELUDIUM project (contract number UMO-2021/41/N/ST1/02908). We appreciate the feedback from two anonymous reviewers, which enhanced the quality of this work.

\addcontentsline{toc}{section}{Acknowledgements}

\bibliographystyle{splncs04}
\bibliography{references.bib}

\end{document}